\documentclass{amsart}
\usepackage{amsmath,amsthm,bbm}
\usepackage{amsfonts,amssymb}

\newtheorem{thm}{Theorem}[section]
\newtheorem{cor}[thm]{Corollary}
\newtheorem{lem}[thm]{Lemma}

\newtheorem{defn}[thm]{Definition}

\theoremstyle{remark}
\newtheorem{rem}{Remark}[section]

\newcommand{\1}{{\mathbbm 1}}

\def\hb{\hfill$\Box$}
\def\f{\frac}

\def\va{\varepsilon}
 \def\a{{\alpha}}
 
 \def\g{{\gamma}}
 \def\k{{\kappa}}
 \def\t{{\theta}}
 \def\l{{\lambda}}
 \def\d{{\delta}}
 \def\o{{\omega}}
 
 \def\la{{\langle}}
 \def\ra{{\rangle}}
 \def\ve{{\varepsilon}}

 \def\CH{{\mathcal H}}

 \def\cC{{\mathcal C}}

 \def\CM{{\mathcal M}}

 \def\CV{{\mathcal V}}

 \def\RR{{\mathbb R}}
 \def\ZZ{{\mathbb Z}}
        
        \def\proj{\operatorname{proj}}
        
\def\bl{\bigl}
\def\br{\bigr}
\def\dmax{\displaystyle \max}

\def\Bl{\Bigl}
\def\Br{\Bigr}
\def\dsup{\displaystyle\sup}
\def\({\Bigl(}
\def\){\Bigr)}
\def\ta{\theta}
\newcommand{\wt}{\widetilde}

\DeclareMathOperator{\meas}{meas}

\begin{document}

\title[maximal function and multiplier theorem on sphere]
{Maximal function and Multiplier Theorem for Weighted Space on the Unit Sphere}
\author{Feng Dai}
\address{Department of Mathematical and Statistical Sciences\\
University of Alberta\\, Edmonton, Alberta T6G 2G1, Canada.}
\email{dfeng@math.ualberta.ca}
\author{Yuan Xu}
\address{Department of Mathematics\\ University of Oregon\\
    Eugene, Oregon 97403-1222.}\email{yuan@math.uoregon.edu}

\thanks{ The first  author  was partially supported  by the
NSERC Canada under grant G121211001. The second author was
partially supported by the National Science Foundation under Grant
DMS-0604056}

\begin{abstract}
For a family of weight functions invariant under a finite reflection group,  the
boundedness of a maximal function on the unit sphere is established and
used to prove a multiplier theorem for the orthogonal expansions with respect
to the weight function on the unit sphere.  Similar results are also established
for the weighted space on the unit ball and on the standard simplex.
\end{abstract}

\maketitle

\section{Introduction}
\setcounter{equation}{0}

The purpose of this paper is to study the maximal function in the
weighted spaces on the unit sphere and the related domains. Let
$S^d =\{x: \|x\| =1\}$ be the unit sphere in $\RR^{d+1}$, where
$\|x\|$ denotes the usual Euclidean norm.  Let $\langle
x,y\rangle$ denote the usual Euclidean inner product. We consider
the weighted space on $S^d$ with respect to the measure
$h_\kappa^2 d \omega$, where $d \omega$ is the surface (Lebesgue)
measure on $S^d$ and the weight function $h_\kappa$ is defined by
\begin{equation}\label{G-weight}
h_\kappa(x) = \prod_{v \in R_+} |\langle x, v\rangle|^{\kappa_v}, \qquad
   x \in \RR^{d+1},
\end{equation}
in which $R_+$ is a fixed positive root system of $\RR^{d+1}$, normalized so
that $\langle v, v \rangle =2$ for all $v \in R_+$, and $\kappa$ is a
nonnegative multiplicity function $v \mapsto \kappa_v$ defined on $R_+$ with
the property that $\kappa_u = \kappa_v$ whenever $\sigma_u$, the reflection
with respect to the hyperplane perpendicular to $u$, is conjugate to
$\sigma_v$ in the reflection group $G$ generated by the reflections
$\{\sigma_v:v \in R_+\}$. The function $h_\kappa$ is invariant under the
reflection group $G$. The simplest example is given by the case
$G=\ZZ_2^{d+1}$ for which $h_\kappa$ is just the product weight function
\begin{equation}  \label{weight}
  h_\k(x) = \prod_{i=1}^{d+1} |x_i|^{\k_i}, \qquad \k_i \ge 0,\
  \    x=( x_1,\cdots, x_{d+1}).
\end{equation}
Denote by $a_\k$  the normalization constant, $a_\k^{-1} =
\int_{S^d} h_\k^2(y) d\omega(y)$.   We consider the weighted space
$L^p(h_\k^2;S^d)$ of functions on $S^d$ with the finite norm
$$
  \|f\|_{\k,p} :=  \Big(a_\kappa \int_{S^d} |f(y)|^p h_\kappa^2(y) d\omega(y)
\Big)^{1/p}, \qquad  1 \le p < \infty,
$$
and for $p = \infty$ we assume that $L^\infty$ is replaced by $C(S^d)$, the
space of continuous functions on $S^d$ with the usual uniform norm
$\| f \|_\infty$.

The weight function \eqref{G-weight} was first studied by Dunkl in the context
of $h$-harmonics, which are orthogonal polynomials with respect to $h_\k^2$.
A homogeneous polynomial is called an $h$-spherical harmonics if it is
orthogonal to all polynomials of lower degree with respect to the inner product
of $L^2(h_\kappa^2; S^d)$. The theory of $h$-harmonics is in many ways
parallel to that of ordinary harmonics (see \cite{DX}). In particular, many
results on the spherical harmonics expansions have been extended to
$h$-harmonics expansions, see \cite{Dai,D2, DX, LX, X05a,X05b} and the
references therein. Much of the analysis of $h$-harmonics depends on the
intertwining operator $V_\k$ that intertwines between Dunkl operators, which
are a  commuting family of first order differential-difference operators, and
the usual partial derivatives.  The operator $V_\k$ is a uniquely determined
positive linear operator.  To see the importance of this operator,
let $\CH_n^{d+1}(h_\k^2)$ denote the space of $h$-harmonics of degree $n$;
the reproducing kernel of $\CH_n^{d+1}(h_\k^2)$ can be written in terms of
$V_\k$ as
\begin{equation}\label{reprod}
 P_n^h(x,y) = \frac{n+\l_k}{\l_\k}
      V_\k \left[ C_n^{\l_k}(\la x , \cdot \ra) \right ](y), \qquad x, y \in S^d,
\end{equation}
where $C_n^\l$ is the $n$-th Gegenbauer polynomial, which is orthogonal
with respect to the weight function $w_\l(t) := (1-t^2)^{\l -1/2}$ on $[-1,1]$,
and
\begin{equation}\label{lambda}
\l_\k =  \g_\k + \frac{d-1}{2} \quad \hbox{with} \quad
     \g_\k = \sum_{v\in R_+} \k_v.
\end{equation}
Furthermore, using $V_\k$, a maximal function that is particularly suitable for
studying the $h$-harmonic expansion is defined in \cite{X05b} by
\begin{equation} \label{eq:CMf}
  \CM_\kappa f(x) :=  \sup_{0< \theta \le \pi} \frac{
  \int_{S^d} |f(y)| V_\kappa \left[\chi_{B(x,\theta)}\right](y)
      h_\kappa^2(y)d\omega(y)}
      {\int_{S^d} V_\kappa \left[\chi_{B(x,\theta)}\right](y)
       h_\kappa^2(y)d\omega(y)},
\end{equation}
where $B(x,\theta): = \{y\in B^{d+1}: \langle x,y\rangle \ge \cos
\theta\}$, $B^{d+1}: = \{x: \|x\|\le 1\} \subset \RR^{d+1}$, and
$\chi_E$ denotes the characteristic function of the set $E$. A
weak type (1,1) inequality was established for $\CM_\k f$ in
\cite{X05b}. The result, however, is weaker than the usual weak
type (1,1) inequality and it does not imply the strong $(p,p)$
inequality. One of our main results in this paper is to establish
a genuine weak type (1,1) result, for which we rely on the general
result of \cite{Stein} on semi-groups of operators. Furthermore,
the Fefferman-Stein type result
$$
 \bigg \|\Bigl(\sum_j |\CM_\k f_j|^2\Bigr)^{1/2}\bigg\|_{\k,p} \le c
  \bigg \|\Bigl(\sum_j |f_j|^2\Bigr)^{1/2} \bigg\|_{\k,p}
$$
also holds, which can be used to derive a multiplier theorem for $h$-harmonic
expansions, following the approach of \cite{BC}. These results are presented
in Section 2.

In the case of $\ZZ_2^{d+1}$, the explicit formula of $V_\k$ as an integral
operator is known, which allows us to link the maximal function $\CM_\k f$
with the weighted Hardy-Littlewood maximal function defined by
\begin{equation} \label{eq:Mf}
   M_k f(x) : = \sup_{0 < \t \le \pi} \frac{ \int_{c(x,\t)} |f(y)| h_\k^2(y)d\omega(y)}
                 { \int_{c(x,\t)}   h_\k^2(y)d\omega(y)},
\end{equation}
where $c(x,\theta) : = \{y \in S^d: \langle x,y\rangle \ge \cos
\theta\}$ is the spherical cap. We will show that the maximal
function $\CM_\k f$ is bounded by a sum of the Hardy-Littlewood
maximal function $M_\k f$. As a consequence, we establish a
weighted weak (1, 1) result for $\CM_kf (x)$, in which the weight
is also of the form \eqref{weight} but with different parameters.
Furthermore, we show that the Fefferman-Stein type inequality
holds in the weighted $L^p$ norm. These results are discussed in
Section 3.

The analysis on the sphere is closely related to the analysis on
the unit ball $B^d$ and on the standard simplex $T^d$. In fact,
much of the results on the later two cases can be deduced from
those on the sphere (see \cite{DX,X05a,X05b} and the references
therein). In particular,  maximal functions are also defined on
$B^d$ and $T^d$ in terms of the generalized translation operators
(\cite{X05b}). We will extend our results on the sphere in Section
2 to these two domains, including a multiplier theorem for the
orthogonal expansions in the weighted space on $B^d$ and $T^d$, in
Sections 4 and 5, respectively.

Throughout this paper, the constant $c$ denotes a generic constant,
which depends only on the values of $d$, $\kappa$ and other fixed parameters
and whose value may be different from line to line. Furthermore, we write
$A\sim B$ if $A \le c B$ and $B\le c A$.

\section{Maximal function  and multiplier theorem on $S^d$}
\setcounter{equation}{0}

\subsection{Background}

In this subsection we give a brief account of what will be needed later
on in the paper. For more background and details, we refer to
\cite{DX,X05a,X05b}.

{\it $h$-harmonic expansion.} Let $\CH_n^{d+1}(h_\k^2)$ denote the space of
spherical $h$-harmonics of degree $n$. It is known that $\dim \CH_n^{d+1}(h_\k^2)
 = \binom{n+d+1}{n} - \binom{n+d-1}{n-2}$. The usual Hilbert space theory shows
that
$$
L^2(h_k^2; S^d) = \sum_{n=0}^\infty \CH_n^{d+1}(h_\k^2): \qquad
       f = \sum_{n=0}^\infty \proj_n^\k f,
$$
where $\proj_n^\k : L^2(h_\k^2;S^d) \mapsto \CH_n^{d+1}(h_\k^2)$ is the
projection operator, which can be written as an integral operator
\begin{equation} \label{projection}
   \proj_n^\k f(x) = a_\k \int_{S^d} f(y) P_n^h(x,y) h_\k^2(y) d\omega(y),
\end{equation}
where $P_n^h$ is the reproducing kernel of $\CH_n^{d+1}(h_\k^2)$, which
satisfies the compact representation \eqref{reprod}.

\medskip
{\it Intertwining operator.}
For a general reflection group, the explicit formula of $V_\k$ is not known.
In the case of $\ZZ_2^{d+1}$, it is an integral operator given by
\begin{equation} \label{Vk}
  V_\kappa f(x) = c_\k
       \int_{[-1,1]^{d+1}} f(x_1 t_1, \ldots,x_{d+1} t_{d+1})
        \prod_{i=1}^{d+1} (1+t_i)(1-t_i^2)^{\kappa_i -1} d t,
\end{equation}
where $c_\k$ is the normalization constant determined by $V_\kappa 1 =1$.
If some $\kappa_i =0$, then the formula holds under the limit relation
$$
 \lim_{\lambda \to 0} c_\lambda \int_{-1}^1 f(t) (1-t)^{\lambda -1} dt
  = [f(1) + f(-1)] /2.
$$

\medskip
{\it Convolution.} For $f \in L^1(h_\k^2;S^d)$ and $g \in
L^1(w_{\l_\k};  [-1,1])$, define (\cite[p.6, Definition
2.1]{X05a})
\begin{equation} \label{convolution}
f \star_\k g(x) :=a_\k \int_{S^d} f(y) V_\k[g(\la x,\cdot \ra)](y)
          h_\k^2(y) d\omega(y).
\end{equation}
This convolution satisfies the usual Young's inequality (see
\cite[p.6, Proposition  2.2]{X05a}): For $f\in L^q(h_\k^2;S^d)$
and $g \in L^r (w_{\l_\k};[-1,1])$, $\|f\star_\k g\|_{k,p} \le
\|f\|_{k,q} \|g\|_{w_{\l_\k},r}$, where $p,q,r \ge 1$ and $p^{-1}
= r^{-1}+q^{-1} -1$. For $\k =0$, $V_\k = id$, this becomes the
classical convolution on the sphere (\cite{CZ}). Notice that by
\eqref{reprod} and \eqref{projection}, we can write $\proj_n^\k f
$ as a convolution.

\medskip

{\it Ces\`aro $(C,\delta)$ means.} For $\delta > 0$, the $(C,\delta)$ means,
$s_n^\delta$, of a sequence $\{c_n\}$ are defined by
$$
s^\delta_n = (A_n^\delta)^{-1}\sum_{k=0}^n A_{n-k}^\delta c_k,
  \qquad A_{n-k}^\delta = \binom{n-k+\delta}{n-k}.
$$
We denote the $n$-th $(C,\delta)$ means of the $h$-harmonic expansion by
$S_n^\delta(h_\kappa^2;f)$. These means can be written as
$$
 S_n^\delta(h_\k^2;f) = (f\star_\k q_n^\delta)(x), \qquad
   q_n^\delta(t) = (A_n^\delta)^{-1}\sum_{k=0}^n A_{n-k}^\delta
      \frac{(k+\l)}{\l} C_k^\l(t),
$$
where $\l = \l_\k$. The function $q_n^\delta(t)$ is the kernel
of the $(C,\delta)$ means of the Gegenbauer expansions at $x =1$.

\medskip

{\it Generalized translation operator $T_\t^\k$}.  This operator
is defined implicitly by  (\cite[p.7]{X05a})
\begin{equation}\label{transl}
c_\l \int_0^\pi T_\t^\k f(x) g(\cos \t) (\sin \t)^{2\l} d\t =  (f\star_\k g)(x),
\end{equation}
where $g$ is any $L^1(w_\l)$ function and $\l = \l_\k$. The operator
$T_\t^\k$ is well-defined and becomes the classical spherical means
$$
    T_\t f(x) = \frac{1}{\sigma_{d} (\sin \theta)^{d-1}}
     \int_{\langle x, y\rangle = \cos\theta} f(y) d\omega(y),
$$
when $\k=0$, where $\sigma_{d} = \int_{S^{d-1}} d\omega = 2
\pi^{d/2}/\Gamma(d/2)$ is the surface area of $S^{d-1}$.
Furthermore, $T_\t^\k$ satisfies similar properties as those
satisfied by $T_\t$, as shown in \cite{X05a,X05b}. In particular,
if $f(x) =1$, then $T_\theta^\kappa f(x) =1$.

\medskip

{\it Spherical caps}. Let $d (x,y) := \arccos \la x ,y \ra$ denote the geodesic
distance of $x,y \in S^d$. For $0 \le \theta  \le \pi$, the set
$$
c(x,\theta) := \{y \in S^d:d(x,y)\le \theta\} =
      \{y \in S^d: \langle x,y\rangle \ge \cos \theta\}
$$
is called the spherical cap centered at $x$.  Sometimes we need to consider
the solid set under the spherical cap, which we denote by $B(x,\theta)$ to
distinguish it from $c(x,\theta)$; that is,
$$
  B(x,\theta): = \{y\in B^{d+1}: \langle x,y\rangle \ge \cos
  \theta\},
$$where $B^{d+1}=\{y\in\mathbb{R}^{d+1}:\  \ \|y\|\leq 1\}$.

\medskip

{\it Maximal function.} For $f\in L^1(h_\kappa^2)$, define (\cite{X05b})
$$
  \CM_\kappa f(x) := \sup_{0< \theta \le \pi}
    \frac{\int_0^\theta T_\phi^\kappa |f|(x) (\sin \phi)^{2\lambda_\k} d\phi}
      {\int_0^\theta (\sin \phi)^{2\lambda_\k} d\phi}.
      $$
      This maximal function can be used to study the $h$-harmonic
expansions, since we can often prove $|(f\star_\k g)(x)| \le c
\CM_\k f(x)$. Using \eqref{transl} it is shown in \cite{X05b} that
an equivalent definition for $\CM_\k f$ is \eqref{eq:CMf}; that
is,
\begin{equation} \label{eq:CMff}
  \CM_\kappa f(x) =  \sup_{0< \theta \le \pi} \frac{
  \int_{S^d} |f(y)| V_\kappa \left[\chi_{B(x,\theta)}\right](y)
      h_\kappa^2(y)d\omega(y)}
      {\int_{S^d} V_\kappa \left[\chi_{B(x,\theta)}\right](y)
       h_\kappa^2(y)d\omega(y)}.
\end{equation}
We note that setting $f(x) =1$ and $g(t) = \chi_{[\cos\t,1]}(t)$ in \eqref{transl}
leads to
\begin{equation}\label{denoM1f}
 a_\kappa \int_{S^d} V_\kappa [\chi_{B(x,\t)}](y)h_\kappa^2(y)d\omega(y)
     = c_{\lambda_\k} \int_0^\t (\sin \phi)^{2\l_\k} d\phi \sim \t^{2\l_\k+1}.
\end{equation}


\subsection{Maximal Function}

To state the weak type inequality, we define, for any measurable subset $E$ of
$S^d$,  the measure with respect to $h_\k^2$ as
$$
     \meas_\k E : = \int_{E} h_\k^2(y) d\o(y).
$$
Our main result in this section is the boundeness of $\CM_\k f$:

\begin{thm} \label{thm: weak(1,1)}
 If $f\in L^1(h_\k^2;S^d)$, then $\CM_\k f$ satisfies
\begin{equation} \label{weak}
   \meas_\k \{x: \CM_\k f(x) \ge \a\}   \le c \frac{\|f\|_{\k
   ,1}}{\a},\  \   \  \forall\a>0.
\end{equation}
Furthermore, if $f \in L^p(h_\k^2;S^d)$ for  $1 < p \le \infty$, then
$\|\CM_k f\|_{\k,p} \le c \|f\|_{\k,p}$.
\end{thm}

The inequality \eqref{weak} is usually refereed to as weak type
(1,1) inequality. In order to prove this theorem, we follow the
approach of \cite{Stein} on general diffusion semi-groups of
operators on a measure space.
For this we need the Poisson integral with respect to $h_\k^2$,
which can be written as \cite[p. 190, Theorem 5.3.3]{DX}
\begin{equation}\label{Poisson}
    P_r^\k f(x) = f \star_\k p_r^\k, \qquad \hbox{where} \quad
           p_r^\k(s) = \frac{1-r^2}{(1-2 r s + r^2)^{\l_k+1}}.
\end{equation}
The kernel $p_r^\k$ is one of the generating function of the Gegenbauer
polynomials of parameter $\l_\k$. Hence, by \eqref{reprod},
we can write $P_r^\k f$ as
$$
 P_r^\k f(x) = \sum_{n=0}^\infty r^n \proj_n^\k f(x), \qquad 0 \le r <1.
$$
from which it follows easily that  $T^t: = P_r^\k f$ with $r =
e^{-t}$ defines a semi-group. Since $V_\k$ is positive and
$p_r^\k$ is clearly nonnegative, $P_r^\k f \ge 0$ if $f \ge 0$. We
see that the semi-group $P_{e^{-t}}^\k f$ is positive. We will
need another semi-group, which is the discrete analog of the heat
operator:
\begin{equation}\label{heat}
H_t^\k f : = f \star_\k q_t^\k, \qquad q_t^\k (s):= \sum_{n=0}^\infty
      e^{-n(n+2\l_\k) t} \frac{n+\l_\k}{\l_\k} \,C_n^{\l_\k}(s).
\end{equation}
In fact,  the $h$-harmonics in $\CH_n^{d+1}(h_\k^2)$ are the
eigenfunctions of an operator $\Delta_{h,0}$, which is the
spherical part of a second order differential-difference operator
analogous to the ordinary Laplacian, the eigenvalues are
$-n(n+2\l_k)$. It follows immediately from $(\ref{heat}$) that
$\{H_t^\k\}_{t\ge 0}$ is a semi-group.
 The following result is the key
for the proof of the theorem.

\begin{lem} \label{lem:P-H}
The Poisson and the heat semi-groups are connected by
\begin{equation}\label{P-H}
  P_{e^{-t}}^\k f (x) =  \int_0^\infty \phi_t(s)  H^\k_s f(x)  ds,
        \end{equation} where $$ \phi_t(s) :=  \frac{t}{2 \sqrt{\pi}}  s^{-3/2}
             e^{-(\frac{t}{2\sqrt{s}} -  \l_\k \sqrt{s})^2}.$$
Furthermore, assume that $f(x) \ge 0$ for all $x$, then for all $t
>0$,
\begin{equation}\label{P-Hmaxima}
 P_*^\k f(x):= \sup_{0<r<1}  P_r^\k f(x) \le
       c  \sup_{s > 0}  \frac{1}{s} \int_0^s H_u^\k f (x) du .
\end{equation}
Consequently, $P_*^\k f$ is bounded on $L^p(h_\k^2;S^d)$ for $1 <
p\le \infty$ and of weak type $(1,1)$.
\end{lem}

\begin{proof} That $\{H_t^\k\}_{t\ge 0}$ is a semi-group is obvious.  Moreover,
since $V_\k$ is positive and $q_t^\k$ is known to be non-negative
(\cite{KM}), it follows that $H_t^\k f$ is positive. The
positivity shows that $\|q_t^\k\|_{w_{\l_\k},1} =1$, so that
$\|H_t^\k f\|_{\k,p} \le \|f\|_{\k,p}$, $1 \le p \le \infty$, by
applying Young's inequality on $f\star_\k q_t^\k$. Thus, using the
Hopf-Dunford-Schwartz ergodic theorem ( \cite[p.48]{Stein}), we
conclude that  the maximal operator $\dsup_{s
> 0} \(\frac{1}{s} \int_0^s H_u^\k f (x) du\)$ is bounded on
$L^p(h_\k^2; S^d)$ for $1<p\leq \infty$ and of week type $(1, 1)$.
Therefore,  it is sufficient   to prove ($\ref{P-H}$) and
($\ref{P-Hmaxima}$).

First we prove ($\ref{P-H}$).  Applying the well known
identity (\cite[p.46]{Stein})
$$
   e^{- v } = \frac{1}{\sqrt{\pi}} \int_0^\infty \frac{e^{-u}} {\sqrt{u}}e^{-v^2/4u} du,
   \qquad v > 0,
$$
with $v = (n+\l_\k) t$ and making a change of variable $s =
t^2/4u$, we conclude that
\begin{align*}
  e^{-nt} &  = e^{ \l_\k t} \frac{1}{\sqrt{\pi}} \int_0^\infty
  \frac{e^{-u}} {\sqrt{u}}
          e^{-\frac{n (n+2\l_\k)t^2}{4u}} e^{-\frac{\l_\k^2t^2}{4u}} du\\
          & =  \frac{t}{2 \sqrt{\pi}} \int_0^\infty  e^{-n (n+2\l_\k)s}s^{-3/2}
             e^{-(\frac{t}{2\sqrt{s}} -  \l_\k \sqrt{s})^2} ds\\
             &=\int_0^\infty e^{-n (n+2\l_\k)s} \phi_t(s)\, ds.
\end{align*}
Multiplying by $\proj_n^\k f$ and summing up over $n$ proves the
integral relation ($\ref{P-H}$).

For the proof of \eqref{P-Hmaxima},  we use ($\ref{P-H}$) and
integration by parts to obtain
\begin{align*}
 P_{e^{-t}}^\k f(x) &\ = - \int_0^\infty  \left( \int_0^s H_u^\k f(x) du \right) \phi'_t(s) ds \\
      &   \le \sup_{s > 0}  \left( \frac{1}{s} \int_0^s H_u^\k f(x)du \right)
             \int_0^\infty s |\phi'_t(s)|ds,
\end{align*}
where the derivative of $\phi'_t(s)$ is taken with respect to $s$.
Also, we  note that by ($\ref{Poisson}$) and
($\ref{convolution}$),
$$\sup_{0<r\leq e^{-1}} P_r^\k f(x)\leq c \|f\|_{1,\k}=c\lim_{s\to\infty} \f 1s \int_0^s H_u^{\k} (|f|)(x)\, du.$$
Therefore, to finish the proof of \eqref{P-Hmaxima}, it suffices to
show that $
\sup_{0<t\leq 1} \int_0^\infty s
|\phi_t'(s)|ds $ is bounded by a constant.

A quick computation shows that $\phi_t'(s)> 0$ if $s < \a_t$ and
$\phi_t'(s) < 0$ if $s > \a_t$, where
$$
\a_t:= \frac{t^2}{3+\sqrt{9+4\l_\k^2t^2}} \sim t^2, \qquad
              \hbox{  $0\leq t\leq 1$}.
$$
Since the integral of $\phi_t(s)$ over $[0,\infty)$ is 1 and $\phi_t(s) \ge 0$,
integration by parts gives
\begin{align*}
  \int_0^\infty s |\phi_t'(s)| ds &\ = 2 \a_t \phi_t(\a_t) - \int_0^{\a_t} \phi_t(s) ds
     + \int_{\a_t}^\infty \phi_t(s)ds \\
       &\    \le 2 \a_t \phi_t(\a_t) +1=\frac{t}{\sqrt{\pi\a_t}
       }e^{-\frac{(t-2\l_\k\a_t)^2}{4\a_t}}+1
        \le c
\end{align*}
 as desired.
\end{proof}

We are now in a position to prove Theorem \ref{thm: weak(1,1)}.

\medskip\noindent
{\it Proof of Theorem \ref{thm: weak(1,1)}.}  From the definition
of $p_r^\k$, it follows easily that if $1 - r \sim \t$, then
\begin{align*}
 p_r^\k(\cos \t) & =
     \frac{1-r^2}{ \left((1- r)^2 + 4 r \sin^2 \frac{\t}{2}\right)^{\l_\k+1}} \\
   & \ge c  \frac{1-r^2}{ \left((1- r)^2 + r \t^2 \right)^{\l_\k+1}}
      \ge c\, (1-r)^{-(2 \l_\k+1)}.
\end{align*}
For $j \ge 0$ define $r_j: = 1-2^{-j}\t$ and set
$
     B_j := \left \{ y \in B^{d+1}: 2^{-j-1}\t \le d(x,y) \le 2^{-j} \t\right \}.
$
The lower bound of $p_r^\k$ proved above shows that
$$
      \chi_{B_j}(y) \le c \, (2^{-j} \t)^{2\l_k+1}p^\k_{r_j}(\la x, y\ra),
$$
which implies immediately that
$$
 \chi_{B(x,\t)}(y) \le \sum_{j=0}^\infty \chi_{B_j}(y)
    \le c \, \t^{2 \l_k+1} \sum_{j=0}^\infty 2^{-j(2\l_\k+1)}p_{r_j}^\k(\la x, y\ra).
$$
Since $V_\k$ is a positive linear operator, applying $V_\k$ to the above
inequality gives
\begin{align*}
  & \int_{S^{d-1}}|f(y)|V_\k \left[ \chi_{B(x,\t)}\right ](y)h_\k^2(y) d\o(y) \\
   & \qquad \le  c \, \t^{2\l_\k+1} \sum_{j=0}^\infty 2^{-j(2\l_\k+1)}
      \int_{S^{d-1}}|f(y)|V_\k \left [ p_{r_j}(\la x, y\ra) \right ](y)h_\k^2(y) d\o(y)\\
   & \qquad   = c \, \t^{2\l_\k+1} \sum_{j=0}^\infty 2^{-j(2\l_\k+1)} P_{r_j}^\k (|f|;x) \\
   & \qquad  \le c  \t^{2\l_\k+1}  \sup_{0< r <1} P_r^\k (|f|; x).
\end{align*}
Dividing by $ \t^{2\l_\k+1} $ and using \eqref{denoM1f}, we have
proved that $\CM_\k f (x) \le c P_*^\k |f|(x)$. The desired result
now follows from Lemma \ref{lem:P-H}.
\qed 

\medskip

A weighted maximal function, call it $M_\k f$, on $\RR^d$ is defined in \cite{TX}
in terms of a translation that is defined via Dunkl transform, the analogue of
Fourier transform for the weighted $L^2(h_\k^2;\RR^d)$. The translation can
be expressed in term of $V_\k$ when acting on radial functions. The boundedness
of the maximal function $M_\k f$ was established in \cite{TX}.  Although the
relation between the maximal function $M_\k f$ and $\CM_\k f$ is not known
at this moment, it should be pointed out that our proof of
Theorem \ref{thm: weak(1,1)} follows the line of argument used  in the
proof of \cite{TX}.


\subsection{A multiplier theorem}

As an application of the above result we state a multiplier theorem.  Let
$\Delta g(t) = g(t+1)  - g(t)$ and $\Delta^k = \Delta^{k-1} \Delta$.

\begin{thm} \label{multiplier}
Let $\{\mu_j\}_{j=0}^\infty$ be a sequence of real  numbers that
satisfies
\begin{enumerate}
\item{}  $\displaystyle\sup_j |\mu_j| \le  c < \infty,$
\item{} $\displaystyle \sup_j 2^{j(k-1)} \sum_{l= 2^j}^{2^{j+1}}|\Delta^k u_l | \le c < \infty$,
\end{enumerate}
where $k$ is the smallest integer $\ge \l_\k+1$.  Then $\{\mu_j\}$
defines an $L^p(h_\k^2;S^d)$, $1<p<\infty$, multiplier; that is,
$$
 \bigg \| \sum_{j=0}^\infty \mu_j \proj_j^\k f \bigg \|_{\k,p} \le c \| f\|_{\k,p}, \qquad 1 < p <\infty,
$$
where $c$ is independent of $\mu_j$ and $f$.
\end{thm}

When $\k =0$, the theorem becomes part  of \cite[Theorem 4.9]{BC}
on the ordinary spherical harmonic expansions. The proof of this
theorem follows that of the theorem in \cite{BC}. One of the main
ingredient is the Littlewood-Paley function
\begin{equation}\label{g(f)}
 g(f)= \left( \int_0^1 (1-r) \left|\frac{\partial}{\partial r}
             P_r^\k f \right|^2 dr \right)^{1/2},
\end{equation}
where $P_r^\k f $ is the Poisson integral with respect to $h_\k^2$ defined
in \eqref{Poisson}. A general Littlewood-Paley theory was established in
\cite{Stein} for a family of diffusion semi-group of operators \{$T^t\}_{t \ge 0}$
on a measure space, in which the $g$ function is defined as
$$
g_1(f) = \left( \int_0^\infty t \left|\frac{\partial}{\partial t} T^t f \right|^2 dt \right)^{1/2}.
$$
Applying the general theory to $T^t: = P_r^\k f$ with $r = e^{-t}$ and using the
fact that the crucial point in the definition of $g(f)$ is when $r$ close to $1$, it
follows that
\begin{equation}\label{Littlewood-Paley}
c^{-1} \|f\|_{\k,p}\leq \|g(f) \|_{\k,p} \le c \|f\|_{\k,p},
\qquad 1 < p <\infty,
\end{equation}
for $f \in L^p(h_\k^2;S^d)$, where the inequality in the left side holds under
the additional assumption that $\int_{S^d} f(y)h_\k^2(y)\,dy=0$. Another
ingredient of the proof is the Ces\`aro means. Recall that the $(C,\delta)$
means are denoted by $S_n^\d(h_\k^2;f)$. What is needed is the following result:

\begin{thm} \label{thm:cesaro}
For $\d> \l_\k$, $1 < p < \infty$ and any sequence  $\{n_j\}$ of positive integers,
\begin{equation} \label{cesaro}
 \Bigg \| \bigg( \sum_{j=0}^\infty
    \big | S_{n_j}^\d (h_\k^2; f_j )\big|^2\bigg)^{1/2} \Bigg \|_{\k,p}
 \le c   \Bigg \| \bigg( \sum_{j=0}^\infty \big| f_{j} \big|^2 \bigg)^{1/2}  \Bigg \|_{\k,p}.
\end{equation}
\end{thm}

\begin{proof}
The proof of \eqref{cesaro} follows the approach of \cite[p.104-5]{Stein} that
uses a generalization of the Riesz convexity theorem for sequences of
functions. Let $L^p(\ell^q)$ denote the space of all sequences $\{f_k\}$ of
functions for which the norm
$$
 \|(f_k)\|_{L^p(\ell^q)}:=\Bigg(\int_{S^d}  \bigg(\sum_{j=0}^\infty |f_j(x)|^q\bigg)^{p/q}
      h_\k^2(x) d\o(x) \Bigg)^{1/p}
$$
are finite. If $T$ is bounded as operator on $L^{p_0}(\ell^{q_0})$
and on $L^{p_1}(\ell^{q_1})$, then the Riesz convexity theorem states
that $T$ is also bounded on $L^{p_t}(\ell^{q_t})$, where
$$
   \frac{1}{p_t} = \frac{1-t}{p_0} + \frac{t}{p_1}, \quad
     \frac{1}{q_t} = \frac{1-t}{q_0} + \frac{t}{q_1}, \quad 0 \le t \le 1.
$$
We apply this theorem on the operator $T$ that maps the sequence
$\{f_j\}$ to the sequence $\{ S_{n_j}^\d(h_\k^2;f_j)\}$. It is
shown in \cite[p.76 and 78]{X05b} that $\sup_n | S_n^\d(h_\k^2;
f)(x)| \le c \CM_\k f(x)$ for all $x \in S^d$ if $\d > \l_\k$.
Consequently, $T$ is bounded on $L^p(\ell^p)$. It is also bounded
on $L^p(\ell^\infty)$ since
$$
   \bigg \|\sup_{j\ge 0} \big|S_{n_j}^\d(h_\k^2;f_j) \big|\bigg\|_{\k,p} \le
     c \bigg \| \CM_\k \Big(\sup_{j\ge 0} |f_j|\Big) \bigg\|_{\k,p}
     \le c \bigg \|\sup_{j\ge 0} |f_j |\bigg\|_{\k,p}.
$$
Hence, the Riesz convexity theorem shows that $T$ is bounded on
$L^p(\ell^q)$ if $1 < p \le q \le \infty$. In particular, $T$ is
bounded on $L^p(\ell^2)$ if $1 < p\leq 2$. The case $2 < p<\infty$
follows by the standard duality argument, since the dual space of
$L^{p}(\ell^2)$ is $L^{p'}(\ell^2)$, where $1/p+1/p' =1$,  under
the paring
$$
  \la (f_j),(g_j)\ra : = \int_{S^d} \sum_j f_j(x)g_j(x) h_\k^2(x) d\o(x)
$$
and $T$ is self-adjoint under this paring as $S_n^\d(h_\k^2)$ is
self-adjoint in $L^2(h_\k^2;S^d)$.
\end{proof}

Using the two ingredients, \eqref{Littlewood-Paley} and
\eqref{cesaro}, the proof of Theorem \ref{multiplier} follows from
the corresponding proof in \cite{BC} almost verbatim.

\begin{rem}
In the case of $\k=0$, the condition $\d > \l_\k=(d-1)/2$ is the critical
index for the convergence of $(C,\d)$ means  in $L^p(h_\k^2; S^d)$
for all $1 \le p \le \infty$. For $h_\k^2$ given in \eqref{weight} and
$G = \ZZ_2^d$, this remains true if at least one $\k_i$ is zero.
However, if $\k_i \ne 0$ for all $i$, then the critical index is
$\l_\k - \min_{1\le i\le d+1}\k_i$ (\cite{LX}). It remains to be seen if the
condition $k \ge   \l_\k +1$ in Theorem \ref{multiplier} can be
improved to $k\ge  \l_\k -\min_{1\leq i\leq d+1} \k_i +1$.
\end{rem}

The proof of Theorem $\ref{thm:cesaro}$ actually yields the
following Fefferman-Stein type inequality (\cite{FS}) for the maximal
function $\mathcal{M}_\k f$.

\begin{cor} \label{cor2-5}
Let $1<p\leq 2$ and $f_j$ be a sequence of functions. Then
\begin{equation}\label{F-S-2-5}
    \bigg \| \Big ( \sum_j (\CM_\k f_j)^2 \Big)^{1/2}\bigg \|_{\k,p} \le
        c \bigg \| \Big(\sum_j |f_j|^2\Big)^{1/2} \bigg \|_{\k,p}.
\end{equation}
\end{cor}
We do not know if the inequality \eqref{F-S-2-5} holds for
$2<p<\infty$  under a general finite reflection group.  However,
it will be shown  in the next section that  \eqref{F-S-2-5} is
true for all $1<p<\infty$  in the case of $G= \ZZ_2^{d+1}$.

\section{Maximal function for product weight}
\setcounter{equation}{0}

The result on the maximal function in the previous section is established for
every finite reflection group. In the case of $G= \ZZ_2^{d+1}$, the weight
function becomes \eqref{weight}; that is,
$$
  h_\k(x) = \prod_{i=1}^{d+1} |x_i|^{\k_i}, \qquad \k_i \ge 0.
$$
We know the explicit formula of the intertwining operator $V_\k$, as shown
in \eqref{Vk}. This additional information turns out to offer more insight into
the maximal function $\CM_k f$. The main result in this section relates
$\CM_\k f$ to the weighted Hardy-Littlewood maximal function.

\begin{defn}
For $f \in L^1(h_\k^2;S^d)$, the weighted Hardy-Littlewood maximal function
is defined by
\begin{equation}\label{HLM}
   M_\k f(x) : = \sup_{0 < \t \le \pi} \frac{ \int_{c(x,\t)} |f(y)| h_\k^2(y)d\omega(y)}
                 { \int_{c(x,\t)}   h_\k^2(y)d\omega(y)}.
\end{equation}
\end{defn}

Since $h_\k$ is a doubling weight (\cite{Dai}),  $M_\k f$ enjoys
the classical properties of maximal functions. We will show that
the maximal function $\CM_\k f$ is bounded by a sum of $M_\k f$,
so that the properties of $\CM_k f$ can be deduced from those of
$M_\k f$. We shall  need several lemmas. The first lemma is an
observation made in \cite[p. 72]{X05b}, which we state as a lemma
to emphasize  its importance in the development below.

\begin{lem} \label{lem1}
For $x \in S^d$ let $\bar x := (|x_1|, \ldots,|x_d|)$. Then the support set of
the function $V_\kappa \left[\chi_{B(x,\theta)}\right](y)$ is
$\{y: d(\bar x, \bar y) \le \theta\}$; in other words,
$$
V_\kappa \left[\chi_{B(x,\theta)}\right](y) = 0   \qquad \hbox{if
   $\langle \bar x, \bar y \rangle < \cos \theta$}.
$$
\end{lem}

\begin{proof}
The explicit formula ($\ref{Vk}$) of $V_\k$ shows that if $V_\k
\left[\chi_{B(x,\theta)}\right](y)
 = 0$ if
$$
\chi_{_{B(x,\t)}}(t_1y_1,t_2y_2, \cdots, t_{d+1} y_{d+1} ) = 0$$
for every
 $t \in [-1,1]^{d+1}$ or    if $x_1y_1t_1 + \ldots + x_{d+1}y_{d+1} t_{d+1} < \cos \t$, which
clearly holds if $\langle \bar x, \bar y \rangle < \cos \theta$.
\end{proof}

Our second lemma contains the essential estimate for an upper bound of
$\CM_\k f$.

\begin{lem} \label{lem2}
For $0 \le \theta \le \pi$, $x=(x_1,\cdots, x_{d+1})\in S^d$ and
$y\in S^d$,
\begin{equation} \label{eq:lem2}
     \left |V_\k \left[ \chi_{B(x,\t)} \right](y)  \right| \le c
           \prod_{j=1}^{d+1} \frac{\t^{2 \k_j}}{(|x_j|+\t)^{2\k_j}} \,
                 \chi_{c(\bar x, \t)}(\bar y).
\end{equation}
\end{lem}

\begin{proof}
The presence of $\chi_{c(\bar x,\t)}(\bar y)$ in the right hand side of the
stated estimate comes from Lemma \ref{lem1}. Hence, we only need to
derive the upper bound of $V_\k\left[\chi_{B(x,\t)} \right](y)$ for
$d (\bar x, \bar y)\le \t$, which we assume in the rest of the proof.
If $\pi/2\le \t \le \pi$, then $\t / (|x_j|+\t) \ge c$
and the inequality \eqref{eq:lem2} is trivial. So we can assume $0 \le \t \le
\pi/2$ below.  By the definition of $V_\k$,
$$
  V_\k\left[\chi_{B(x,\t)} \right](y) = c_\k \int_{\sum_{i=1}^{d+1}  t_ix_iy_i \ge \cos\t}
       \prod_{i=1}^{d+1} (1+t_i) (1-t_i^2)^{\k_i-1}dt
$$
where $t$ also satisfies  $t\in [-1,1]^{d+1}$. We first enlarge
the domain of integration to $\{t \in [-1,1]^{d+1}:
\sum_{i=1}^{d+1}|t_ix_iy_i| \ge \cos \t\}$ and replace $(1+t_i)$
by 2, so that we can use the $\ZZ_2^{d+1}$ symmetry of the
resulted integrant to replace the integral to the one on
$[0,1]^{d+1}$,
\begin{align*}
  V_\k\left[\chi_{B(x,\t)} \right](y) \le \ & c  \int_{\sum_{i=1}^{d+1} |t_ix_iy_i| \ge \cos\t}
       \prod_{i=1}^{d+1} (1-t_i^2)^{\k_i-1}dt \\
       \le \ & c  \int_{t\in[0,1]^{d+1},  \sum_{i=1}^{d+1} t_i |x_iy_i| \ge \cos\t}
                   \prod_{i=1}^{d+1} (1-t_i^2)^{\k_i-1}dt.
\end{align*}
To continue, we enlarge the domain of the integral to $\{t \in [0,1]^{d+1}:
 t_j|x_jy_j| +\sum_{i\ne j} |x_iy_i| \ge \cos \t\}$ for each fixed $j$ to obtain
$$
  V_\k\left[\chi_{B(x,\t)} \right](y) \le c \prod_{j=1}^{d+1}
      \int_{0 \le t_j \le 1, t_j|x_jy_j| +\sum_{i\ne j} |x_iy_i| \ge \cos \t }
          (1-t_j)^{\k_j-1}dt_j.
$$
For each $j$ we denote the last integral by $I_j$.  First of all, there is the
trivial estimate
$
       I_j \le  \int_0^1 (1-t_j)^{\k_j-1} dt_j = \k_j^{-1}.
$
Secondly, for $\la \bar x, \bar y\ra \ge \cos \t$, we have the estimate
$$
  I_j \le \int_{\frac{\cos\t - \sum_{i\ne j} |x_iy_i|}{|x_jy_j|}}^1 (1-t_j)^{\k_j-1} dt_j
     = \k_j^{-1} \frac{(\la \bar x, \bar y\ra - \cos \t)^{\k_j}}{|x_jy_j|^{\k_j}}.
$$
Together, we have established the estimate
$$
  I_j \le \k_j^{-1} \min
   \left\{1, \frac{(\la \bar x, \bar y\ra - \cos \t)^{\k_j}}{|x_jy_j|^{\k_j}} \right\}.
$$
Recall that $d (\bar x, \bar y)\le \t$. Assume first that $|x_j| \ge 2 \t$. Then
$|x_j| \ge (|x_j|+\t)/2$. The inequality $\big | |x_j| - |y_j| \big | \le
\|\bar x - \bar y\| \le d(\bar x,\bar y) \le \t$ implies that $|y_j| \ge |x_j| - \t
\ge |x_j|/2$, so that $|y_j| \ge (|x_j|+\t)/4$. Furthermore, write $t :=
d(\bar x, \bar y) \le \t$ and recall that $\t \le \pi/2$. We have then
$$
\la \bar x, \bar y \ra - \cos \t = \cos t  - \cos \t = 2
      \sin \tfrac{\t - t}{2} \sin \tfrac{t+\t}{2} \le (\t-t) \t \le \t^2.
$$
Putting these ingredients together, we arrive at an upper bound for $I_j$,
$$
 I_j \le c \frac{\t^{2\k_j}} {(|x_j|+\t)^{2\k_j}},
$$
under the assumption that $|x_j| \ge 2 \t$. This estimate also holds for
$|x_j| \le 2 \t$, since in that case $\t / (|x_j|+\t) \ge 1/3$. Thus, the last
inequality holds for all $x$ and for all $j$, from which the stated inequality
follows immediately.
\end{proof}

Our next lemma gives the order of the denominator in $M_\k f$,
which was proved in \cite[(5.3), p. 157]{Dai} in the case when
$\displaystyle\min_{1\leq j \leq d+1} \tau_j \ge 0$.

\begin{lem}\label{lem3}
Let $\tau=(\tau_1,\cdots, \tau_{d+1})\in\mathbb{R}^{d+1}$ with
$\displaystyle \min_{1\leq j \leq d+1} \tau_j >-1$. Then for $0
\le \t \le \pi$ and $x=(x_1,\cdots, x_{d+1})\in S^d$,
$$
   \int_{c(x,\t)}  \prod_{j=1}^{d+1} |y_j|^{\tau_j}  d\o(y) \sim
     \t^d \prod_{j=1}^{d+1} (|x_j|  +\t)^{\tau_j},
$$where the constant of equivalence depends only on $d$ and
$\tau$.
\end{lem}

\begin{proof}
Without loss of generality we may assume that  $x_j\ge 0$ for all
$1 \le j \le d+1$ and $x_{d+1}=\dmax_{1\leq j \leq d+1} x_j$, as well
as $0<\t < \frac1{2\sqrt{d+1}}$. Since $x_{d+1} ={ \dmax_{1\leq j \leq d+1} }x_j \ge
\frac 1{\sqrt{d+1}}$,  it follows that
\begin{equation}\label{3-2}
     y_{d+1}\ge x_{d+1}-\t \ge
     \frac 1{2\sqrt{d+1}},\  \  \  \forall y=(y_1, \cdots, y_{d+1})\in
      c (x,\t).
\end{equation}
Using \eqref{3-2} and the fact that $d\o(y)= c_d (1-\|\wt{y}\|^2)^{-\f12}
\, d\wt{y}$ for  $y=(\wt{y}, y_{d+1})$ and $y_{d+1}=\sqrt{1-\|\wt{y}\|^2}\ge  0$,
as well as the fact that $|x_j-y_j| \le \|x-y\|\le d(x,y)$, we conclude
\begin{align*}
  \int_{c(x,\t)}  \prod_{j=1}^{d+1}
       |y_j|^{\tau_j}  d\o(y) & \sim  \int_{d(x,y)\leq \t}
         \prod_{j=1}^{d} |y_j|^{\tau_j} dy_1\, dy_2\, \cdots dy_d\\
   &  \leq c \prod_{j=1}^d \int_{x_j-\ta} ^{x_j+\ta} |y_j|^{\tau_j} \,
        dy_j \sim  \t^d \prod_{j=1}^{d+1} (|x_j|
   +\t)^{\tau_j}.
\end{align*}
where in the last step, if $\tau_j < 0$, consider the cases $x_j \ge 2 \t_j$ and
$x_j < 2 \t_j$ separately. This gives the desired upper estimate.

For the proof of the lower estimate, let $z = (z_1,\ldots,z_{d+1}) \in S^d$ be
defined by $z_j=x_j+\va\t$ for $j=1,2,\cdots, d$ and  $ z_{d+1} =
(1-z_1^2-\cdots-z_d^2)^{\f12}$, where $\va>0$ is a sufficiently small constant
depending only on $d$. Using \eqref{3-2}, a quick computation shows that
$$
  \|x-z\|^2 = d (\va \t)^2 + \frac{|z_{d+1}^2- x^2_{d+1}|^2}{(z_{d+1}+ x_{d+1})^2}
     \le  d(\va \t)^2 + (d+1)d (2  \va \t +  \va^2\t^2)^2,
$$
from which and the fact that $2 \sin \frac{d(x,z)}{2} = \|x-z\|$, it follows that we
can choose $\va$ small enough so that $z \in c(x, \frac \t2)$.   Consequently,
$c(z, \frac{\va \t}2)\subset c(x,\ta)$ and, for any  $y=(y_1,\cdots, y_{d+1})\in
c(z, \frac{\va\t}2)$,
$$
x_j + \frac { \va \t}2 = z_j -\frac { \va \t}2\leq y_j\leq z_j
+ \frac { \va \t}2 =x_j + \frac { 3\va \t}2,\   \  \ j=1,2,\cdots,
d+1,
$$
which implies immediately that
$$
\prod_{j=1}^{d+1} |y_j|^{\tau_j} \sim \prod_{j=1}^{d+1}
    (|x_j|+\t)^{\tau_j}, \qquad  \forall y\in c(z, \frac{\va\t}2),
$$
and, as a consequence,
$$
 \int_{c(x,\t)}  \prod_{j=1}^{d+1} |y_j|^{\tau_j}  d\o(y)\ge
    \int_{c(z,\frac{\t\va}2)}  \prod_{j=1}^{d+1} |y_j|^{\tau_j}  d\o(y)
       \ge c \, \t^d \prod_{j=1}^{d+1} (|x_j| +\t)^{\tau_j},
$$
proving  the desired lower estimate.
\end{proof}

In particular, Lemma $\ref{lem3}$ shows that $h_\k^2$ is a
doubling weight in the sense that
$$
  \meas_\k c(x, 2\t) \le c \meas_\k c(x,\t), \qquad \forall x \in S^d, \quad \t >0.
$$

We are now ready to prove our first main result. For $x\in \RR^{d+1}$ and
$\ve \in \ZZ_2^{d+1}$, we write $x \ve: = (x_1 \ve_1,\ldots,x_{d+1}\ve_{d+1})$.

\begin{thm} \label{Mf}
Let $f \in L^1(h_\k^2;S^d)$. Then for every  $x \in S^d$,
\begin{equation} \label{MfCMf}
     \CM_\k f(x) \le c \sum_{\ve\in \ZZ_2^{d+1}} M_\k f(x\ve).
\end{equation}
\end{thm}

\begin{proof}Since
$$
   \bl\{y\in S^d:\  \  d(\bar{x}, \bar{y})\leq \ta\br\}
   =\bigcup_{\varepsilon\in\ZZ_2^{d+1}} \bl\{ y\in S^d:\  \
      d(x\varepsilon, y)\leq \t\br\},
$$
it follows from Lemmas \ref{lem1} that
\begin{align*}
   J_\t f(x) &\ : =  \int_{S^d} |f(y)| V_\k\left[\chi_{B(x,\t)}\right](y) h_\k^2(y)d\o(y) \\
     &\  =  \int_{\la  \bar x, \bar y \ra \ge \cos\t} |f(y)|
         V_\k\left[\chi_{B(x,\t)}\right](y) h_\k^2(y)d\o(y)\\
          &\ \le  \sum_{\ve \in \ZZ_2^{d+1}}
    \int_{\la x\ve , y\ra\ge \cos\t} |f(y)| V_\k\left[\chi_{B(x,\t)}\right](y) h_\k^2(y)d\o(y).
\end{align*}
Consequently, using Lemmas \ref{lem2} and \ref{lem3}, we conclude that
\begin{align*}
  J_\t f(x) &\  \le c \sum_{\ve \in \ZZ_2^{d+1}}
       \prod_{j=1}^{d+1}\frac{\t^{2\k_j}} {(|x_j|+\t)^{2\k_j}}
           \int_{\la x\ve,y\ra\ge \cos\t} |f(y)|  h_\k^2(y)d\o(y)\\
   &\ \le c \  \t^{2|\k|+ d} \sum_{\ve \in\ZZ_2^{d+1}} M_\k f(x\ve).
 \end{align*}
Dividing the above inequality by $\t^{2 |\k|+d} =\t^{2\l_\k+1}$ and, recall
\eqref{denoM1f}, taking the supremum over $\t$ lead to \eqref{MfCMf}.
\end{proof}

There are several applications of Theorem \ref{Mf}. First we need several
notations. For $x=(x_1,\cdots, x_{d+1})$, $y=(y_1,\cdots, y_{d+1})\in
\mathbb{R}^{d+1}$, we write $x<y$  if $x_j<y_j$ for all $1\leq j\leq d+1$.
We denote by $\1$  the vector $\1:= (1,1,\cdots,1)\in\mathbb{R}^{d+1}$.
Moreover, we extend  the definitions of  $h_\tau$, $\meas_\tau$,
$\|\cdot\|_{\tau,p}$, $L^p(h_\tau^2;S^d)$ and  $M_\tau$  to the full range
of $\tau=(\tau_1,\cdots,\tau_{d+1})>-\f {\1}2$.  Thus,
$$
h_\tau(x)=\prod_{j=1}^{d+1} |x_j|^{\tau_j},\qquad
    \|f\|_{\tau, p} =\left(\int_{S^d} |f(x)|^p h_\tau^2(x)\, d\o(x)\right)^{\f1p}
$$
and $M_\tau$ denotes the Hardy-Littlewood maximal function with
respect to the measure $h_\tau^2(x)\, d\o(x)$, as defined in
($\ref{HLM}$).

As an  application of Theorem \ref{Mf}, we can  prove  the
boundedness of $\CM_\k f$ on  the  spaces  $L^p(h_\tau^2; S^d)$
for a wider range of $\tau$ without using  the
Hopf-Dunford-Schwartz ergodic theorem.

\begin{thm} \label{thm: MCM}
 If  $-\f{\1}2< \tau \leq \k$ and $f\in
L^1(h_\tau^2;S^d)$, then $\CM_\k f$ satisfies
\begin{equation}\label{3-4-2}
   \meas_\tau \{x: \CM_\k f(x) \ge \a\}   \le c
   \frac{\|f\|_{\tau,1}}{\a},\   \    \   \forall\a>0.
\end{equation}
Furthermore, if  $1 < p \leq \infty$, $-\f{\1} 2 <\tau<p\k +
\f{p-1}2\1$ and $f \in L^p(h_\tau^2;S^d)$,  then
\begin{equation}\label{3-5-2}\|\CM_k f\|_{\tau,p} \le c
\|f\|_{\tau,p}.\end{equation}
\end{thm}

\begin{proof}
We start with the proof of ($\ref{3-4-2}$).  Note that  if
$\tau=(\tau_1,\cdots,\tau_{d+1})\leq \k$, then
$$
 \int_{c(x,\ta)} |f(y)|h_\k^2(y)\, d\o(y) \leq c
  \bigg(\prod_{j=1}^{d+1} (|x_j|+\t)^{2(\k_j-\tau_j)}\bigg)\int_{c(x,\ta)}
  |f(y)|h_\tau^2(y)\, d\o(y),
$$
which, together with Lemma $\ref{lem3}$, implies
$$
M_\k f(x) \leq c M_\tau f (x),\   \   \  x\in S^d,\   \ \tau\leq \k.
$$
Hence, using  the inequality \eqref{MfCMf}, we
obtain that ,  for $-\f{\1} 2< \tau \leq \k$,
\begin{align*}
   \meas_\tau\{x: \CM_k f(x) \ge \a\} &\ \le \sum_{\ve \in \ZZ_2^{d+1}}
       \meas_\tau \left \{x:  M_\k f(x \ve) \ge  c \, \a/ 2^{d+1}\right \} \\
&\le \sum_{\ve \in \ZZ_2^{d+1}}
       \meas_\tau \left \{x:  M_\tau f(x \ve) \ge  c' \, \a\right \} \\
&\ = \sum_{\ve \in \ZZ_2^{d+1}}
       \int_{\left \{y: \  M_\tau f(y\ve) \ge  c' \, \a \right\} }
       h_\tau^2(y) d\o(y) \\
     &\ = 2^{d+1} \int_{\left \{x :  M_\tau f(x) \ge  c' \, \a\right \} }
      h_\tau^2(y) d\o(y)\\
     &\ \le  c  \frac{\|f\|_{\tau,1} } {\a}
\end{align*}
where we have used the  $\mathbb{Z}_2^{d+1}$- invariance of
$h_\tau$ in the fourth step, and the fact that $M_\tau $ is of
weak type (1,1) with respect to the doubling measure
$h_\tau^2(y)\, d\o(y)$ in the last step. This proves
($\ref{3-4-2}$).

For the proof of  ($\ref{3-5-2}$), we choose a number  $q\in (1,p)$
such that $\tau< q \k +\f{q-1}2\1$ and claim that it is sufficient to
prove
\begin{equation}\label{3-6-1}
     M_\k f(x) \leq c \( M_\tau(|f|^q)(x)\)^{\f 1q}.
 \end{equation}
Indeed, using \eqref{3-6-1}, the inequality ($\ref{3-5-2}$) will follow
from \eqref{MfCMf}, the $\mathbb{Z}_2^{d+1}$ invariance of $h_\tau$
and the boundedness of the maximal function $M_\tau$ on the space
$L^{\f pq}(h_\tau^2; S^d)$.

To prove ($\ref{3-6-1}$), we use H\"older's inequality with $q'=\f q{q-1}$
and Lemma $\ref{lem3}$ to obtain
\begin{align*}
\int_{c(x,\t)} &|f(y)|h^2_\k (y)\, d\o(y) \\
    & \leq \left(\int_{c(x,\t)}
      |f(y)|^q h_\tau^2(y) \, d\o(y) \right)^{\f1q} \left(\int_{c(x,\t)}
     h^2_{q'\k-\f{q'}q\tau} (y)\, d\o(y)\right)^{\f1{q'}}\\
   &\sim \left(\int_{c(x,\t)} |f(y)|^q h_\tau^2(y) \, d\o(y)\right)^{\f1q}
     \Bigg(\prod_{j=1}^{d+1}
    (|x_j|+\t)^{2\k_j-\f{2\tau_j}q}\Bigg)\t^{\f{d}{q'}}\\
   &\sim \meas_\k (c(x,\t))\left( \f
    1{\meas_{\tau}(c(x,\t))}\int_{c(x,\t)} |f(y)|^q h_\tau ^2(y) \, dy
   \right)^{\f1q},
\end{align*}
where we have also used the fact that the assumption $\tau< q \k +\f{q-1}2\1$
is equivalent to $q'\kappa - \frac{q'}{q} \tau > - \frac{\1}{2}$. This
proves ($\ref{3-6-1}$)  and completes the proof.
\end{proof}

For our next application of Theorem \ref{Mf} we will need the following result.

\begin{lem} \label{stein_lem}
Let $1 < p < \infty$ and let $W$ be a nonnegative, local integrable function
on $S^d$. Then
\begin{equation}\label{3-10-1}
\int_{S^d} | M_\k f(x)|^p W(x) h_\k^2(x)\, d\o(x)\leq c_p
     \int_{S^d}|f(x)|^p M_\k W(x)h_\k^2(x)\, d\o(x).
\end{equation}
\end{lem}

Such a result was first proved in \cite{FS} for maximal function on $\RR^d$.
The proof can be adopted easily to yield Lemma \ref{stein_lem}. Indeed, the
fact that $h_\k^2$ is a doubling weight shows that the Hardy-Littlewood
maximal function defined by \eqref{HLM} satisfies
$$
   M_\k f(x) \sim \sup_{x \in E \in \cC} \frac{\int_{E} |f(y)| h_\k^2(y)d\o(y)}
                              {\int_{E}  h_\k^2(y)d\o(y)},
$$
where $\cC$ is the collection of all spherical caps in $S^d$, which implies that
$$
  \int_{c(x,\t)} |f(y)|h_\k^2(y) dy \le c \left( \meas_\k c(x,\t)\right)
             \inf_{z\in c(x,\t)} M_\k f(z)
$$
for any spherical cap $c(x,\t)$. As a consequence, we can prove the key
inequality
$$
   \meas_\k (E) \le  \frac{c}{\a} \int_{S^d} |f(y)| M_\k W(y) h_\k^2(y) d\o(y)
$$
for any compact set $E$ in $\{x \in S^d: M_k f(x) > \alpha\}$, as in the proof
for the maximal function on $\RR^d$ in \cite[p. 54-55]{Stein95}.  In fact,
\eqref{3-10-1} holds with $h_\k^2 (y)d\o$
replaced by any doubling measure $d\mu$ on the sphere.

An important tool in harmonic analysis is the Fefferman-Stein type inequality
(\cite{FS}), which we established in Corollary \ref{cor2-5} for $\CM_\k f$
in the case of $1<p \le 2$ and a general reflection group. In the current
setting of $G=\ZZ_2^{d+1}$, we can use Theorem \ref{Mf} to prove a weighted
version of this inequality for $1< p< \infty$.

\begin{thm} \label{cor2}
Let $1<p<\infty$, $-\f{\1}2<\tau <p\k + \f
{p-1}2\1$,  and let $\{f_j\}_{j=1}^\infty$ be a sequence
of functions. Then
\begin{equation}\label{F-S}
    \Bigg \| \bigg( \sum_{j=1}^\infty (\CM_\k f_j)^2 \bigg)^{1/2}\Bigg \|_{\tau,p} \le
        c \Bigg \| \bigg(\sum_{j=1}^\infty |f_j|^2\bigg)^{1/2} \Bigg \|_{\tau,p}.
\end{equation}
\end{thm}

\begin{proof}
Using Theorem \ref{Mf} and the Minkowski inequality, we obtain
\begin{align*}
   & \Bigg\| \bigg ( \sum_j \left(\CM_\k f_j\right)^2 \bigg)^{1/2}  \Bigg \|_{\tau,p}
        \ \le c  \Bigg \| \bigg( \sum_j  \bigg( \sum_{\ve \in \ZZ_2^{d+1}}
                                M_\k f_j (x\ve) \bigg)^2 \bigg)^{1/2} \Bigg \|_{\tau,p} \\
    & \qquad \ \le c  \sum_{\ve \in \ZZ_2^{d+1}}  \Bigg \|
               \(\sum_j \left(M_\k f_j (x\ve) \right)^2 \)^{1/2} \Bigg \|_{\tau,p}
            \ \le c  \Bigg \|  \(\sum_j (M_\k f_j)^2\)^{1/2} \Bigg \|_{\tau,p}.
\end{align*}
Thus, it is sufficient to prove
\begin{equation}
\label{F-S-1}
     \Bigg \|  \bigg( \sum_j (M_\k f_j)^2 \bigg)^{1/2} \Bigg \|_{\tau,p} \le
        c  \Bigg \| \bigg( \sum_j |f_j|^2\bigg)^{1/2}  \Bigg \|_{\tau,p}.
\end{equation}

We start with the case $1<p\leq 2$. Let $q$ be chosen such that
$1<q<p$ and $\tau<q\k+\f{(q-1)\1}2$. We use the inequality
($\ref{3-6-1}$) to obtain
\begin{align*}
 \Bigg \| \bigg ( \sum_j \left(M_\k f_j\right)^2 \bigg)^{1/2}\Bigg \|_{\tau,p}
  &\ \leq c \Bigg \| \bigg(\sum_j \left( M_\tau
     (|f_j|^q)\right)^{\f2q}\bigg)^{\f q2}\Bigg\|_{\tau, \f pq}^{\f1q}\\
 & \leq c \Bigg\|\bigg(\sum_j |f_j|^{q\cdot \f 2q}\bigg)^{\f q2}
 \Bigg\|_{\tau, \f pq}^{\f1q}
 =\Bigg \|\bigg(\sum_j |f_j|^2\bigg)^{\f12}\Bigg\|_{\tau, p},
\end{align*}
where we have used the classical Fefferman-Stein inequality for
the maximal function $M_\tau$  and the space $L^{\f pq}(\ell^{\f
2q})$ in the second step. This proves ($\ref{F-S-1}$) for $1<p\leq
2$.

Next, we consider the case $2<p<\infty$. Noticing  that
$$
-\f{\1}2<\tau<p\k+\f {(p-1)\1}2\Longleftrightarrow
-\f{\1}2<\f2p\tau+ \left(\f 1p-\f12 \right)\1 < 2\k +\f{\1}2,
$$
we may choose a vector $\mu\in\mathbb{R}^{d+1}$ such that
\begin{equation}\label{3-9-1}
  -\f{\1}2<\f2p\tau+\left(\f 1p-\f12\right)\1<\mu<2\k +\f{\1}2,
\end{equation}
and a number $1<q<2$ such that $\mu< q\k+\f{q-1}2\1$.  Let $g$
be a nonnegative function on $S^d$  satisfying $\|g\|_{\tau, \f p{p-2}}=1$
and
$$
\Bigg\|\bigg(\sum_j |M_\k f_j|^2\bigg)^{\f12}\Bigg\|_{\tau, p}^2 =\int_{S^d}
\bigg(\sum_j |M_\k f_j(x)|^2\bigg) g(x) h_\tau^2(x)\, d\o(x).
$$
Then by the assumption $\mu< q\k+\f{q-1}2\1$, \eqref{3-6-1},
\eqref{3-10-1} with $p = 2/q>1$ and H\"older's inequality, we obtain
\begin{align*}
   \int_{S^d}  \bigg(\sum_j |M_\k f_j(x)|^2\bigg) & g(x) h_\tau^2(x)\, d\o(x)
     \leq c \sum_j \int_{S^d} \(M_\mu \bl( |f_j|^q\br)(x)\)^{\f 2q}
   g(x) h_\tau^2(x) \,  d\o(x)\\
& \leq c \int_{S^d} \bigg(\sum_j |f_j(x)|^2 \bigg)M_\mu( g
   h_{\tau-\mu}^2)(x)h_\mu^2(x)\, d\o(x)\\
 &\leq c \Bigg\|\bigg(\sum_j
   |f_j|^2\bigg)^{\f12}\Bigg\|_{\tau, p}^2
   \Bl\| M_\mu ( g h_{\tau-\mu}^2) h_{\mu-\tau}^2\Br\|_{ \tau, \f{p}{p-2}}.
\end{align*}
Using the boundedness of $M_\k f$ and ($\ref{3-9-1}$), we have
\begin{align*}
\Bl\| M_\mu ( g h_{\tau-\mu}^2) h_{\mu-\tau}^2\Br\|_{ \tau, \f{p}{p-2}}
 & = \Bl\|M_\mu (gh_{\tau-\mu}^2)\Br\|_{ \f{p}{p-2} \mu -\f 2{p-2}
\tau, \f p{p-2}}\\
& \leq c \Bl\|gh_{\tau-\mu}^2\Br\|_{\f{p}{p-2} \mu -\f
2{p-2}\tau, \f p{p-2}}\\
&=c  \|g\|_{\tau, \f{p}{p-2}} = c.
\end{align*}
Putting these two inequalities together, we have proved the inequality
($\ref{F-S-1}$) for  the case $2<p<\infty$.
\end{proof}

\begin{rem}
It is shown in \cite{X05b} that, for $\d > \l_\k$,  the Ces\`aro $(C,\d)$ means
satisfy $|S_n^\d(h_\k^2; f)| \le c \CM_\k f(x)$. Hence, we can get a weighted
inequality for the Ces\`aro means by replacing $\CM_\k f_j$ in \eqref{F-S}
by $S_{n_j}^\d(h_\k^2; f_j)$. This gives a $\|\cdot\|_{\tau,p}$ weighted version of
Theorem \ref{thm:cesaro} that holds under the condition $-\f{\1}2<\tau <p\k + \f
{p-1}2\1$.
\end{rem}


\section{Maximal function and multiplier theorem on $B^d$}
\setcounter{equation}{0}

Analysis in weighted spaces on the unit ball
$B^d=\{x\in\mathbb{R}^d:\   \  \|x\|\leq 1\}$ in $\RR^d$ can often
be deduced from the corresponding results on $S^d$; see
\cite{DX,X05a,X05b} and the reference therein. Below we develop
results analogous to those in the previous sections.

\subsection{Weight function invariant under a reflection group}

Let $\k =(\k', k_{d+1})$ with $\k'=(\k_1, \cdots, \k_d)$ and assume
$\k_i\ge 0$  for $1 \le i \le d+1$. Let $h_{\k'}$ be the weight function
\eqref{G-weight}, but defined on $\RR^d$, that is invariant under a
reflection group $G$. We consider the weight functions on
$B^d$ defined by
\begin{equation}\label{weightB}
  W_{\k}^B(x) := h_{\k'}^2(x) (1-\|x\|^2)^{\k_{d+1}-1/2},    \qquad x \in B^d,
\end{equation}
which is invariant under the reflection group $G$. Under the mapping
\begin{equation}\label{B-S}
\phi: x\in B^d \mapsto (x,\sqrt{1-\|x\|^2}) \in S^{d}_+:=\{y \in S^d:y_{d+1} \ge 0\}
\end{equation}
and multiplying the Jacobian of this change of variables, the weight function
$W_\k^B$ comes exactly from $h_{\k}^2$ defined by
$$
  h_{\k}(x_1,\ldots, x_{d+1}):= h_{\k'}^2(x_1,\ldots,x_d) |x_{d+1}|^{2\k_{d+1}}.
$$
The weight function $h_{\k}$ is invariant under the reflection
group $G \times \ZZ_2$. All of the results established in Section
2 holds for $h_{\k}$.

We denote the $L^p(W_{\k}^B;B^d)$ norm by $\|f\|_{W_\k^B,p}$.
The norm of $g$ on $B^d$ and its extension on $S^d$ are related by the identity
\begin{equation}\label{BSintegral}
\int_{S^d} g(y) d\omega = \int_{B^d} \left[ g(x,\sqrt{1-\|x\|^2}\,)+
     g(x,-\sqrt{1-\|x\|^2}\,) \right]\frac{dx}{\sqrt{1-\|x\|^2}}.
\end{equation}
The orthogonal structure is preserved under the mapping
\eqref{B-S} and the study of orthogonal expansions for $W_{\k}^B$
can be essentially reduced to that of $h_\k^2$. In fact, let
$\CV_n^d(W_{\k}^B)$ denote the space of orthogonal polynomials of
degree $n$ with respect to $W_{\k}^B$ on $B^d$. The orthogonal
projection, $\proj_n(W_{\k}^B; f)$,  of $f \in L^2(W_{\k}^B;B^d)$
onto $\CV_n^d(W_{\k}^B)$ can be expressed in terms of the
orthogonal projection of $F(x,x_{d+1}):= f(x)$ onto
$\CH_n^{d+1}(h_\k^2)$:
\begin{equation}\label{projBS}
    \proj_n(W_{\k}^B; f, x) = \proj_n^\k F(X), \qquad X := (x,\sqrt{1-\|x\|^2}).
\end{equation}
Furthermore, a maximal function was defined in \cite[p.81]{X05b}
in terms of the generalized translation operator of the orthogonal
expansion. More precisely,  let
$$
e(x,\theta): = \{(y,y_{d+1}) \in B^{d+1}:
     \langle x, y\rangle + \sqrt{1-\|x\|^2}\, y_{d+1}
          \ge \cos \theta, \quad y_{d+1} \ge 0\}.
$$
Then this maximal function, denoted by $\CM_{\k}^B f(x)$, was
shown to satisfy the relation
\begin{equation*}
 \CM_{\kappa}^B f(x) = \sup_{0< \theta \le \pi}
   \frac{\int_{B^d}|f(y)| V_{\k}^B \left[\chi_{e(x,\theta)}\right](Y)
      W_{\k}^B(y) dy}
        {\int_{B^d} V_{\k }^B \left[\chi_{e(x,\theta)}\right](Y)
      W_{\k}^B(y) dy},
\end{equation*}
where  $Y= (y,\sqrt{1-\|y\|^2})$, and for $g: \RR^{d+1} \mapsto
\RR$,
$$
  V_{\k}^B g(x,x_{d+1}): =   \frac{1}{2} \left[V_{\k} g(x,x_{d+1})
           + V_{\k} g(x, -x_{d+1}) \right],
$$
in which $V_{\k}$ is the intertwining operator associate with $h_{\k}$.
This maximal function can be written in terms of the maximal function
$\CM_{\k} f$ in \eqref{eq:CMff}. In fact, we have
$\CM_{\k}^B f(x) = \CM_{\k} F(X)$.
Our main result in this section states that $\CM_{\k}^B f$ is of weak (1,1).
Let us define
$$
   \meas_{\k}^B E : = \int_{E} W_{\k}^B(x)dx, \qquad E\subset B^d.
$$

\begin{thm} \label{thm: MCMB}
 If $f \in L^1(W_{\k}^B;B^d)$ then  $\CM_k^B$ satisfies
$$
  \meas_{\k}^B \left\{x \in B^d: \CM_\k^B f(x) \ge \a \right\}
        \le c \frac{\|f\|_{W_{\k}^B,1}}{\a},\   \   \forall  \a>0.
$$
Furthermore, if $f \in L^p (W_{\k}^B;B^d)$ for $1 <p\le \infty$,
then $\|\CM_{\k}^B f \|_{W_{\k}^B,p} \le c \|f\|_{W_{\k}^B,p}$.
\end{thm}

\begin{proof}
Since $\CM_{\k}^B f(x) = \CM_{\k} F(X)$, it follows from
\eqref{BSintegral} that
\begin{align*}
 \meas_{\k}^B \left\{x \in B^d: \CM_{\k}^B f(x) \ge \a \right\} & =
      \int_{B^d} \chi_{\{\CM_{\k}^B f(x) \ge \a \}} (x) W_{\k}^B(x) dx \\
    & = \int_{S_+^d} \chi_{\left\{\CM_{\k}  F(y) \ge \a \right\}} (y)
           h_{\k}^2(y) d\o(y).
\end{align*}
Enlarging the domain of the last integral to the entire $S^d$ shows that
$$
    \meas_{\k}^B \left\{x \in B^d: \CM_\k^B f(x) \ge \a \right\}
     \le  \meas_{\k} \left\{y \in S^{d}: \CM_{\k} F(y) \ge \a \right\}.
$$
Consequently, by Theorem \ref{thm: weak(1,1)}, we obtain
$$
   \meas_{\k}^B \left\{x \in B^d: \CM_{\k}^B f(x) \ge \a \right\}  \le
      c \frac{\|F\|_{\k,1}}{\a}
$$
from which the weak (1,1) inequality follows from $\|F\|_{\k,1} =
\|f\|_{W_{\k}^B,1}$. Since $\CM_{\k}^B f$ is evidently of strong
type $(\infty, \infty)$, this completes the proof.
\end{proof}

The connection \eqref{projBS} and \eqref{BSintegral} allow us to deduce a
multiplier theorem for orthogonal expansion with respect to $W_{\k}^B$
from Theorem \ref{multiplier}.

\begin{thm} \label{multiplierB}
Let $\{\mu_j\}_{j=0}^\infty$ be a sequence  of real numbers that
satisfies
\begin{enumerate}
\item{}  $\dsup_j |\mu_j| \le  c < \infty,$
\item{} $\dsup_j 2^{j(k-1)} \sum_{l= 2^j}^{2^{j+1}}|\Delta^k u_l | \le c < \infty$,
\end{enumerate}
where $k$ is  the smallest  integer $\ge  \l_\k  +1$, and $\l_\k
=\f{d-1}2+\sum_{j=1}^{d+1} \k_j$.  Then $\{\mu_j\}$ defines an
$L^p(W_\k^B;B^d)$, $1<p<\infty$, multiplier; that is,
$$
 \Bigg \| \sum_{j=0}^\infty \mu_j \proj_j^\k f \Bigg \|_{W_\k^B,p} \le c \| f\|_{W_\k^B,p},
      \qquad 1 < p < \infty,
$$
where $c$ is independent of $\{\mu_j\}$ and $f$.
\end{thm}


\subsection{Weight function invariant under $\ZZ_2^d$}
In the case of $G = \ZZ_2^{d}$, the weight function becomes
\begin{equation}\label{weightBprod}
  W_\k^B(x) := \prod_{i=1}^d |x_i|^{2 \k_i} (1-\|x\|^2)^{\k_{d+1}-1/2},
          \qquad x \in B^d,
\end{equation}
which corresponds to the product weight function $h_\k^2(x)=
\prod_{i=1}^{d+1}|x_i|^{2 \k_i}$.
 Taking into the consideration of the
boundary, an appropriate  distance on $B^d$ is defined by
$$
  d_B(x,y) = \arccos \(\la x ,y \ra + \sqrt{1-\|x\|^2}
  \sqrt{1-\|y\|^2}\),\    \    x, y\in B^d,
$$
which is just the projection of the geodesic distance of $S_+^d$
on $B^d$. Thus, one can define the weighted Hardy-Littlewood
maximal function as
$$
   M_\k^B f(x) : = \sup_{0 < \t\le \pi} \frac{\int_{d_B(x,y)\le \t} |f(y)| W_\k^B(y)dy}
         {\int_{d_B(x,y)\le \t} W_\k^B(y)dy}, \   \   x\in B^d.
$$

We have the following analogue of Theorem \ref{Mf}.

\begin{thm} \label{prop3.1}
Let $f\in L^1(W_\k^B;B^d)$. Then for any $x \in B^d$,
\begin{equation}\label{CMBfMF}
   \CM_\k^B f(x) \le c \sum_{\ve \in \ZZ_2^d} M^B_\k f(x \ve).
\end{equation}
\end{thm}

The proof of Theorem \ref{prop3.1} replies on the following lemma,
which implies, in particular, that  $W_\k^B(y)$ is a doubling
weight on $B^d$.

\begin{lem}\label{4-4-2}
If $\tau=(\tau_1,\cdots, \tau_{d+1})>-\f 12\1$, then for
any $x=(x_1,\cdots, x_d)\in B^d$ and $0\leq\t\leq \pi$,
$$ \int_{d_B(y, x)\leq \t} W_\tau^B(y)\, dy \sim
\t^d\prod_{j=1}^{d+1} (|x_j|+\t)^{2\tau_j},$$ where $x_{d+1}
=\sqrt{1-\|x\|^2}$ and $W_\tau^B(y)$ is defined as in
\eqref{weightBprod}.
\end{lem}

\begin{proof}
Recall that $X=(x, x_{d+1})$, and  $c(X,\t)=\{ z\in S^d:\  \ d(X,z)\leq\t\}$.
Set
$$
c_{+}(X,\t)=\left\{(y_1,\cdots, y_{d+1})\in c(X,\t):\  \  y_{d+1}\ge 0\right\}.
$$
From ($\ref{BSintegral}$) it follows that
\begin{equation}\label{4-7-1}
\int_{d_B(y,x)\leq \t} W_\tau^B(y)\, dy
=\int_{c_{+}(X,\t)}h_\tau^2(z)\, d\o(z),
\end{equation}
which, together with Lemma $\ref{lem3}$, implies the desired  upper
estimate
$$
\int_{d_B(y,x)\leq \t} W_\tau^B(y)\, dy\leq \int_{c(X,\t)}h_\tau^2(z)\,
d\o(z)\leq c \, \t^d\prod_{j=1}^{d+1} (|x_j|+\t)^{2\tau_j}.
$$

To prove the lower estimate, we choose a point $z=(z_1,\cdots,
z_{d+1})\in c(X, \f \t2)$ with $z_{d+1}\ge \va \t$, where $\va>0$
is a sufficiently small constant depending only on $d$. Clearly,
$c(z,\f{\va\t}2) \subset c_{+} (X,\t)$.  Hence, by $(\ref{4-7-1})$,
we obtain
\begin{align*}
\int_{d_B(y,x)\leq \t} W_\tau^B(y)\, dy &\ge
\int_{c(z,\f{\va\t}2)} h_\tau^2(y)\, d\o(y) \\
& \sim
\t^d\prod_{j=1}^{d+1} (|z_j|+\t)^{2\tau_j}
\sim \t^d\prod_{j=1}^{d+1} (|x_j|+\t)^{2\tau_j},\end{align*}
where we have used   Lemma $\ref{lem3}$  in the second step, and
the fact that $z\in c(X,\t)$ in the last step.   This gives the
desired lower estimate.
\end{proof}

Now we are in a position to prove Theorem \ref{prop3.1}.\\

\medskip\noindent
{\it Proof of Theorem \ref{prop3.1}.}  It is shown in \cite[p.81]{X05b} that
$$
\int_{B^d} V_\k [ \chi_{e(x,\t)}](Y)W_\k^B (y)\, dy \sim
\t^{2\l_\k +1},
$$ where
$\l_\k=\f{d-1}2+\sum_{j=1}^{d+1} \k_j$.
The proof follows almost exactly as in the proof of Theorem \ref{Mf},
the main effort lies in the proof of the following inequality:
\begin{equation}\label{4-8}
  V_\k^B [\chi_{e(x,\t)}](Y)\leq c \prod_{j=1}^{d+1} \f
   {\t^{2\k_j}}{(|x_j|+\t)^{2\k_j}}\chi_{{\{ y\in B^d:\ d_B (\bar x,\  \bar y)
      \leq \t\}}} (y),
\end{equation}
where
$x_{d+1}=\sqrt{1-\|x\|^2}$ and $\bar z=(|z_1|,\cdots, |z_d|)$ for
$z=(z_1,\cdots, z_d)\in B^d$. However, using ($\ref{Vk}$) and
the fact that $y_{d+1} = \sqrt{1-\|y\|^2}$, we have
\begin{align*}
V_\k^B [\chi_{e(x,\t)}](Y)&=\f12\( V_\k [\chi_{e(x,\t)}](y,
y_{d+1})+V_\k [\chi_{e(x,\t)}](y, -y_{d+1})\)\\
&= c_\k\int_{D}\prod_{j=1}^d (1+t_j)(1-t_j^2)^{\k_j-1}
(1-t_{d+1}^2)^{\k_{d+1}-1}\, dt,
\end{align*}
where
$$
D=\Bl\{ (t_1,\cdots,t_d,  t_{d+1})\in [-1,1]^d \times [0,1]:\   \
\sum_{j=1}^{d+1} t_j x_j y_j \ge \cos\t\Br\}.
$$
This last integral can be estimated exactly as in the proof of Lemma \ref{lem2},
which yields the desired inequality $(\ref{4-8}$).
 \hb

\medskip

As a consequence of Theorem \ref{prop3.1}, we have the following
analogues of Theorems \ref{thm: MCM} and \ref{cor2}.

\begin{cor} \label{thm: MCM-4}
 If  $-\f{\1} 2< \tau \leq \k$ and $f\in
L^1(W^B_\tau;B^d)$, then $\CM_\k f$ satisfies
\begin{equation*}
   \meas^B_\tau \{x: \CM^B_\k f(x) \ge \a\}   \le c
   \frac{\|f\|_{_{W_\tau^B,1}}}{\a},\   \    \   \forall\a>0.
\end{equation*}
Furthermore, if  $1 < p < \infty$, $-\f{\1} 2 <\tau<p\k +
\f{p-1}2\1$ and $f \in L^p(W^B_\tau;B^d)$,  then
\begin{equation*}
\left\|\CM^B_k f\right\|_{W_\tau^B,p} \le c
\|f\|_{W^B_\tau,p}.
\end{equation*}
\end{cor}

\begin{cor} \label{cor2-4}
Let $1<p<\infty$, $-\f{\1}2<\tau <p\k + \f
{p-1}2\1$,  and let $\{f_j\}_{j=1}^\infty$ be a sequence
of functions. Then
\begin{equation*}
    \Bigg \| \bigg( \sum_{j=1}^\infty (\CM^B_\k f_j)^2 \bigg)^{1/2}
        \Bigg \|_{W_\tau^B,p} \le
        c \Bigg \| \bigg(\sum_{j=1}^\infty |f_j|^2\bigg)^{1/2} \Bigg \|_{W^B_\tau,p}.
\end{equation*}
\end{cor}

Using  the formula $\CM_{\k}^B f(x) = \CM_{\k}
F(X)$ and  the method of  \cite{X05b},  one can also deduce
Corollaries \ref{thm: MCM-4} and \ref{cor2-4} directly from
Theorems \ref{thm: MCM} and \ref{cor2}.


\section{Maximal function and multiplier theorem on $T^d$}
\setcounter{equation}{0}

Just like the connection between the structure of function spaces
on $S^d$ and $B^d$, analysis in weighted spaces on the simplex
$$
T^d=\bl\{(x_1,\cdots, x_d)\in\mathbb{R}^d:\  x_1 \ge 0, \ldots, x_d \ge 0 \ \
 \text{and}\  \  x_1+\cdots+x_d\leq 1\br\}
$$ can
often be deduced from the corresponding results on $B^d$; see
\cite{DX,X05a,X05b} and the reference therein.

\subsection{Weight function associated with a reflection group}
Let  $\k'=(\k_1,\cdots, \k_d)$ and   $h_{\k'}$ be the weight
function \eqref{G-weight} on $\RR^d$ invariant under the
reflection group $G$.  We further require that $h_{\k'}$ is even
in all of its variables; in other words, we require that $h_{\k'}$
is invariant under the semi-direct product of $G$ and $\ZZ_2^d$.
Let $\k_{d+1}\ge 0$ and $\k=(\k', \k_{d+1})$.  The weight
functions on $T^d$  we consider  are
\begin{equation} \label{weightT}
  W_\k^T(x) := h_{\k'}(\sqrt{x_1},\ldots,\sqrt{x_d})(1-|x|)^{\k_{d+1}-1/2},
          \qquad x \in T^d,
\end{equation}
where $|x| = x_1+ \ldots + x_d$.  These weight functions are
related to $W_\k^B$ in \eqref{weightB}. In fact, $W_\k^T$ is
exactly the weight function $W_\k^B$ under the mapping
\begin{equation}\label{psi}
\psi: (x_1,\ldots,x_d) \in T^d \mapsto (x_1^2, \ldots, x_d^2) \in B^d
\end{equation}
and upon multiplying the Jacobian of this change of variables. We
denote the norm of $L^p(W_\k^T;T^d)$ by $\|\cdot\|_{W_\k^T,p}$.
The norm of $g$ on $T^d$ and $g\circ \psi$ on $B^d$ are related by
\begin{equation} \label{T-B}
\int_{B^d} g(x_1^2,\ldots,x_d^2) dx = \int_{T^d}
g(x_1,\ldots,x_d)
            \frac{dx}{\sqrt{x_1\cdots x_d}}.
\end{equation}
The orthogonal structure is preserved under the mapping \eqref{psi}. Let
$\CV_n^d(W_\k^T)$ denote the space of orthogonal polynomials of degree
$n$ with respect to $W_\k^T$ on $T^d$. Then $R \in \CV_n^d(W_\k^T)$
if and only if $R\circ \psi \in \CV_{2n}^d(W_\k^B)$. The orthogonal projection,
$\proj_n(W_\k^T; f)$,  of $f \in L^2(W_\k^T;T^d)$ onto $\CV_n^d(W_\k^T)$
can be expressed in terms of the orthogonal projection of $f \circ \psi$ onto
$\CV_{2n}(W_\k^B)$:
\begin{equation}\label{projTB}
   \left ( \proj_n(W_\k^T; f) \circ \psi \right)(x) =
       \frac{1}{2^d} \sum_{\ve \in \ZZ_2^d} \proj_{2n}(W_\k^B; f\circ \psi, x\ve).
\end{equation}
The fact that $\proj_n(W_\k^T)$ of degree $n$ is related to $\proj_{2n}(W_\k^B)$
of  degree $2n$ suggests that some properties of the orthogonal expansions
on $B^d$ cannot be transformed directly to those on $T^d$.

A maximal function $\CM_k^T f$ is defined in \cite[p. 86,
Definition 4.5]{X05b} in terms of the generalized translation
operator of the orthogonal expansion. It is closely related to the
maximal function $\CM_k^B f$ on $B^d$. It was shown in
\cite[Proposition 4.6]{X05b} that
\begin{equation} \label{MfT-B}
          (\CM_\k^Tf )\circ \psi = \CM_\k^B (f\circ \psi).
\end{equation}
We show that this maximal function is of weak (1,1). Let us define
$$
   \meas_\k^T E : = \int_{E} W_\k^T(x)dx, \qquad E\subset T^d.
$$

\begin{thm} \label{thm: MCMT}
If $f \in L^1(W_\k^T;T^d)$. Then  $\CM_k^T$ satisfies
$$
\meas_\k^T \left\{x \in T^d: \CM_\k^T f(x) \ge \a \right\}
            \le c \frac{\|f\|_{W_\k^T,1}}{\a},\   \   \
            \forall\a>0.
$$
Furthermore,  if $f \in L^p(W_{\k}^T;T^d)$ for $1 <p\le \infty$,
then $\|\CM_{\k}^T f \|_{W_{\k}^T,p} \le c \|f\|_{W_{\k}^T,p}$.
\end{thm}

\begin{proof}
Using the relation \eqref{MfT-B} and \eqref{T-B}, we obtain
\begin{align*}
 \int_{T^d} \chi_{\{x \in T^d: \CM_\k^T f(x) \ge \a \}} (x) W_\k^T(x) dx
=  \int_{B^d} \chi_{\{x \in B^d: \CM_\k^T (f\circ \psi)(x) \ge \a  \}} (x) W_\k^B(x) dx.
\end{align*}
Hence, by Theorem \ref{thm: MCMB}, we conclude that
$$
 \meas_\k^B \left\{x \in B^d: \CM_\k^B (f\circ \psi)(x) \ge \a \right\}
  \le c \frac{\|f\circ \psi\|_{W_\k^B,1}}{\a} =  c \frac{\|f\|_{W_\k^T,1}}{\a},
$$
where the last step follows again from \eqref{T-B}.
\end{proof}

The relation \eqref{projTB} shows that we cannot expect to deduce
all results on orthogonal expansion with respect to $W_\k^T$ on
$T^d$ from those on $B_d$.  This applies to the multiplier
theorem. On the other hand, as it is shown in \cite[p.85]{X05b},
we can introduce a convolution $\star_\k^T$ structure and write
$\proj_n^\k (W_\k^T; f) = f \star_\k^T P_n$. Moreover,  we often
have the inequality $|f \star_\k^T g(x)| \le c \CM_\k^T (x)$. For
example, for the Ces\`aro $(C,\d)$ means $S_n^\d (W_\k^T; f)$, we
have
$$
 \sup_n  |S_n^\d (W_\k^T; f,x)| \le c \CM_\k^T (x), \qquad\hbox{if} \quad
             \d > \l_\k=\sum_{j=1}^{d+1}\k_j+\frac{d-1}{2}.
$$
Using this result, we can prove an analogue of Theroem
\ref{thm:cesaro} almost verbatim.  Furthermore, the Poisson
operator,  $P_r^T f$, of the orthogonal expansion with respect to
$W_\k^T$ on  $T^d$ is still a semi-group when we define $T^t f =
P_r^T f$ with $r = e^{-t}$ ( see, for example, \cite[p.90]{X05b}).
So,  the Littlewood-Paley function $g(f)$, defined as in
\eqref{g(f)}, is bounded in $L^p(W_\k^T; T^d)$ for $1 < p <
\infty$. Hence, all the essential ingredients of the proof of the
multiplier theorem in \cite{BC} hold for the orthogonal expansion
with respect to $W_\k^T$. As a consequence, we have the following
multiplier theorem.

\begin{thm} \label{multiplierT}
Let $\{\mu_j\}_{j=0}^\infty$ be a sequence that satisfies
\begin{enumerate}
\item{}  $\dsup_j |\mu_j| \le  c < \infty,$
\item{} $\dsup_j 2^{j(k-1)} \sum_{l= 2^j}^{2^{j+1}}|\Delta^k u_l | \le c < \infty$,
\end{enumerate}
where $k$ is the smallest integer $\ge \l_\k+1$.  Then $\{\mu_j\}$
defines an $L^p(W_\k^T;T^d)$, $1<p<\infty$, multiplier; that is,
$$
 \Bigg \| \sum_{j=0}^\infty  \mu_j \proj_j^\k f \Bigg \|_{W_\k^T,p} \le c \| f\|_{W_\k^T,p},
      \qquad 1 < p <\infty,
$$
where $c$ is independent of $f$ and $\mu_j$.
\end{thm}

\subsection{Weight function associated with $\ZZ_2^d$}
In the case $G = \ZZ_2^d$, we are dealing with the classical weight
function on $T^d$,
\begin{equation} \label{weightTprod}
  W_\k^T(x) := \prod_{i=1}^d |x_i|^{\k_i-1/2} (1-|x|)^{\k_{d+1}-1/2},
          \qquad x \in T^d.
\end{equation}
Under the mapping \eqref{psi}, this weight function corresponds to
$W_\k^B$ at \eqref{weightBprod}. Taking into the consideration of the
boundary, an appropriate distance on $T^d$ is defined by
$$
    d_T(x,y) = \arccos \left ( \la x^{\frac12},y^\frac12 \ra +
         \sqrt{1-|x|} \sqrt{1-|y|} \right), \quad x,y \in T^d,
$$
where $x^{\frac12} = (x_1^\frac12, \ldots,x_d^\frac12)$ for $x \in T^d$.
Evidently, we have $d_B (x,y) = d_T(\psi(x),\psi(y))$. Using this distance,
one can define the weighted Hardy-Littlewood maximal function as
$$
   M_\k^T f(x) : = \sup_{0 < \t\le \pi} \frac{\int_{d_T(x,y)\le \t} |f(y)| W_\k^T(y)dy}
         {\int_{d_T(x,y)\le \t} W_\k^T(y)dy}, \   \   x\in T^d.
$$

We have the following analogue of Theorem \ref{prop3.1}.

\begin{thm} \label{prop4.1}
Let $f\in L^1(W_\k^T;T^d)$. Then for any $x \in T^d$,
\begin{equation}\label{5.7}
        \CM_\k^T f(x) \le c  M^T_\k f(x).
\end{equation}
\end{thm}

\begin{proof}
Using \eqref{T-B}, it follows readily from the definitions of
$M_\k^Bf$ and $M_\k^T f$  that $(M_\k^T f)\circ \psi   = M_\k^B(f
\circ \psi)$. Hence, using the fact that if $g$ is invariant under
the sign changes, then $M_\k^B g (x\va) = M_\k^B g(x)$ by a simple
change of variables, it follows from \eqref{MfT-B} and Theorem
\ref{prop3.1} that
\begin{align*}
   (\CM_\k^T f)\circ \psi (x) & = \CM_\k^B (f \circ \psi) (x)
       \le c \sum_{\va \in \ZZ_2^d} M_\k^B (f\circ \psi) (x \va) \\
       & = c' M_\k^B (f\circ \psi) (x ) = c' (M_\k^T f)\circ \psi (x)
\end{align*}
for $x \in T^d$, from which the stated result follows immediately.
\end{proof}

Although the proof of this theorem may look like a trivial
consequence of the definition of $M_\k f$, we should mention that
the definition of $\CM_\k^T f$ in \cite[p. 86, Definition
4.5]{X05b} is given in terms of the general translation operator
of the orthogonal expansions with respect to $W_\k^T$.

As a consequence of Theorem \ref{prop4.1} or by \eqref{MfT-B} and
\eqref{T-B}, we have the following analogues of Corollaries
\ref{thm: MCM-4} and \ref{cor2-4}:

\begin{cor} \label{thm: MCM-5}
 If  $-\f{\1} 2< \tau \leq \k$ and $f\in
L^1(W^T_\tau;T^d)$, then $\CM_\k f$ satisfies
\begin{equation*}
   \meas^T_\tau \{x: \CM^T_\k f(x) \ge \a\}   \le c
   \frac{\|f\|_{_{W_\tau^T,1}}}{\a},\   \    \   \forall\a>0.
\end{equation*}
Furthermore, if  $1 < p < \infty$, $-\f{\1} 2 <\tau<p\k +
\f{p-1}2\1$ and $f \in L^p(W^T_\tau;T^d)$,  then
\begin{equation*}
\left\|\CM^T_k f\right\|_{W_\tau^T,p} \le c \|f\|_{W^T_\tau,p}.
\end{equation*}
\end{cor}

\begin{cor} \label{cor2-4-5}
Let $1<p<\infty$, $-\f{\1}2<\tau <p\k + \f
{p-1}2\1$,  and let $\{f_j\}_{j=1}^\infty$ be a sequence
of functions. Then
\begin{equation*}
    \Bigg \| \bigg( \sum_{j=1}^\infty (\CM^T_\k f_j)^2 \bigg)^{1/2}
        \Bigg \|_{W_\tau^T,p} \le
        c \Bigg \| \bigg(\sum_{j=1}^\infty |f_j|^2\bigg)^{1/2} \Bigg \|_{W^T_\tau,p}.
\end{equation*}
\end{cor}

\enddocument